\documentclass[11pt]{amsart}
\usepackage{amsmath}
\usepackage{amsthm}
\usepackage{amssymb}
\usepackage{amsfonts,mathrsfs}
\usepackage{stmaryrd}
\usepackage{amsxtra}  
\usepackage{epsfig}
\usepackage{verbatim}

\theoremstyle{plain}
\newtheorem{theorem}[equation]{Theorem}
\newtheorem{proposition}[equation]{Proposition}
\newtheorem{lemma}[equation]{Lemma}
\newtheorem{corollary}[equation]{Corollary}
\newtheorem{definition}[equation]{Definition}
\theoremstyle{remark}
\newtheorem{remark}[equation]{Remark}
\numberwithin{equation}{section}

\newcommand{\cb}{{\mathcal B}}

\newcommand{\cf}{{\mathcal F}}

\newcommand{\co}{{\mathcal O}}

\newcommand{\cv}{{\mathcal V}}

\newcommand{\sa}{{\mathscr A}}
\newcommand{\sB}{{\mathscr B}}

\newcommand{\sg}{{\mathscr G}}

\newcommand{\fR}{{\mathfrak R}}

\newcommand{\C}{{\mathbb C}}

\newcommand{\N}{{\mathbb N}}

\newcommand{\R}{{\mathbb R}}

\newcommand{\Z}{{\mathbb Z}}

\begin{document}

\title[Twisted projections and extension]{Non-holomorphic projections and extension of biholomorphic mappings}
\author{Jeffery D. McNeal}
\thanks{Research partially supported by an NSF grant}
\address{Department of Mathematics, \newline Ohio State University, Columbus, Ohio, USA}
\email{mcneal@math.ohio-state.edu}
\maketitle
\section{{\bf Introduction}}

A fundamental problem in complex analysis is the following: given bounded domains
$\Omega_1, \Omega_2\subset\C^n$ with smooth boundaries and a biholomorphic map $F:\Omega_1\longrightarrow
\Omega_2$, determine conditions on $\Omega_1, \Omega_2$ which guarantee that $F$ and $F^{-1}$
extend smoothly to the closures of the domains, $\overline{\Omega}_1, \overline{\Omega}_2$. When $n=1$, such maps always extend, without further conditions on $\Omega_1$ and $\Omega_2$. In several variables, however, it is unknown whether such maps universally extend or whether there are obstructions to extension. Furthermore, the extension problem takes on additional significance in several variables because the Riemann mapping theorem does not hold: 
when $n>1$, the moduli space of biholomorphism classes of domains is infinite dimensional, even for the subclass of simply connected domains. Positive results about smooth extension to the boundary allows examination of a given equivalence class of domains by studying differential invariants on the boundary of the domains.

A groundbreaking result on this problem was obtained by Fefferman, \cite{fefferman74}, who showed that $F, F^{-1}$ extend smoothly to $\overline{\Omega}_1, \overline{\Omega}_2$ if both domains are strongly pseudoconvex. Fefferman's remarkable proof involved delicate estimates of the Bergman kernel, obtained by analyzing multiple error terms arising from locally approximating the boundaries of the domains, $b\Omega_1, b\Omega_2$ by euclidean balls, and used strong pseudoconvexity in several essential ways.

A subsequent, highly successful approach to this problem was initiated by Bell and Ligocka, \cite{BellLigocka}, and further developed by Bell, \cite{Bell81a}. The Bell-Ligocka program focused on a regularity property of the Bergman projection, rather than strong pseudoconvexity, and eventually led to showing the extension property holds on broad classes of weakly pseudoconvex domains. Let $B=B_\Omega$ denote the Bergman projection on $\Omega$, the orthogonal projection of $L^2(\Omega)$ onto its subspace of holomorphic functions. Say that $\Omega$ satisfies Condition R
if $B:C^\infty(\overline\Omega)\longrightarrow C^\infty(\overline\Omega)$. The main result in \cite{Bell81a} is the following: if $\Omega_1$ satisfies Condition R and $\Omega_2$ is pseudoconvex (both domains having smooth boundary), then $F$ and $F^{-1}$ extend smoothly to $\overline\Omega_1$ and $\overline\Omega_2$, respectively.

The question then arises: which smoothly bounded domains satisfy Condition R? There are many  hypotheses on $\Omega$ known to imply this condition, see \cite{Barrett86}, \cite{BellBoas81}, \cite{BoaStr91}, \cite{Catlin84b}, \cite{Catlin87}. For pseudoconvex domains, these results infer Condition R from a global regularity property of the $\bar\partial$-Neumann operator $N$, specifically that (*)
$N:\Lambda^{0,1}(\overline\Omega)\longrightarrow\Lambda^{0,1}(\overline\Omega)$, where
$\Lambda^{0,1}(\overline\Omega)$ denotes the $(0,1)$ forms with components in $C^\infty(\overline\Omega)$. The $\bar\partial$-Neumann operator is basic operator in complex analysis (see \cite{FolKoh72}) that inverts the $\bar\partial$-Laplacian with natural boundary conditions. The operator $\bar\partial^* N$ gives the special solution to the Cauchy-Riemann equations that is orthogonal to holomorphic functions on $\Omega$ --- this property establishes a relationship between $N$ and $B$. It was an open question, for many years, whether $N$ always satisfied (*) on a smoothly bounded, pseudoconvex domain. However, Christ  gave a negative answer to this question in  \cite{Christ96}, showing that $N$ does not satisfy property (*) on some pseudoconvex domains . Christ's theorem thus limits the applicability of Bell's theorem, though we emphasis that  \cite{Christ96} does not give a counterexample to the smooth extension of biholomorphic mappings.
The earlier works of Barrett, \cite{Barrett92}, and Kiselman \cite{Kiselman} were important precursors to the results in \cite{Christ96}.

The Bergman projection bears on the extension problem through its transformation formula:
\begin{equation}\label{E:bergtransform}
B_1\big(JF\cdot g\circ F\big)=JF\cdot\big[B_2(g)\big]\circ F,\qquad g\in L^2\left(\Omega_2\right)
\end{equation}
where $B_j$ denotes the Bergman projection on $\Omega_j$, $j=1,2$ and $JF$ is the determinant of the holomorphic Jacobian of $F$.
The crucial element in the Bell-Ligocka approach is the fact that $B_2$ has a large null space that is connected to a space of functions reproduced by $B_2$.
Consider an equivalence relation on functions $f,g\in L^2\left(\Omega_2\right)$ defined $f\sim g$ if $B_2f=B_2g$. Let $[f]$ denote the equivalence class of $f$. Then each class $[f]$ contains representatives that vanish to high order on $b\Omega_2$: if $g\in C^\infty\left(\overline{\Omega}_2\right)$ and $M\in\Z^+$ , there exists $g_M\in [g]$ such that $g_M=0$ on $b\Omega_2$ to order $M$.\footnote{This fact, which we refer to as Bell's Lemma, allows functions $g\in A^2\left(\Omega_2\right)$ on the right-hand side of \eqref{E:bergtransform} to be replaced by functions $g_M$ on the left-hand side of  \eqref{E:bergtransform}. But $JF\cdot g_M\circ F$ is essentially smooth up to $b\Omega_1$ if $M$ is large, by Cauchy's estimates. This lets one use \eqref{E:bergtransform} to gain control of $F$'s boundary behavior through Condition R.} The proof of Bell's Lemma hinges on two facts:

\begin{itemize}
\item[(a)] $B_2$ reproduces every function in $A^\infty\left(\Omega_2\right)=\co\left(\Omega_2\right)\cap C^\infty\left(\overline{\Omega}_2\right)$ (as $A^\infty\left(\Omega_2\right)\subset A^2\left(\Omega_2\right)$ since $\Omega_2$ is bounded).
\item[(b)] $B_2$ annihilates a purely anti-holomorphic derivative of every function in $C^\infty\left(\overline{\Omega}_2\right)$ that also vanishes on $b\Omega_2$.
\end{itemize}
It follows that if $f\in C^\infty\left(\overline{\Omega}_2\right)$ is given, functions $\eta_\ell$ of the form $\frac{\partial}{\partial z_k}(r_2\cdot\sigma_\ell)$ with $\sigma_\ell\in C^\infty\left(\overline{\Omega}_2\right)$ can be chosen which have the same Taylor coefficients as $f$, up to order $\ell$, in the variable $r_2$ near $b\Omega_2$. Because of (b), it follows that $B_2(f-\eta_\ell)=B_2(f)$.

There are other solution operators to the Cauchy-Riemann equations besides $\bar\partial^* N$ and some of them are know to have good global regularity properties. Also, some of these operators are connected to projection operators.
Kohn \cite{Kohn73} produced such a solution operator on any smoothly bounded, pseudoconvex domain $\Omega$.
For $t >0$, let
$L^2_t(\Omega)$ be the Hilbert space of functions on $\Omega$ with inner product
$$\langle f,g\rangle_t=\int_\Omega f(z)\overline{g(z)} e^{-t|z|^2}\, dV(z).$$
Kohn showed that for any $s\in\Z^+$, there exists $t_0$ such that the weighted $\bar\partial$-Neumann operator $N_t$ maps $W^s(\Omega)\longrightarrow W^s(\Omega)$ boundedly, if $t >t_0$, where $W^s(\Omega)$ denotes the ordinary Sobolev norm of order $s$.
If $B_{t|z|^2}$ denotes the orthogonal projection of $L^2_t(\Omega)$ onto $\co (\Omega)$ in this inner product, then it follows that
$$B_{t|z|^2}:W^s(\Omega)\longrightarrow W^s(\Omega),\qquad\text{if } t>t_0, $$
boundedly. However, the weighted Bergman projection $B_{t|z|^2}$ does not exhibit the correct connection between the functions it reproduces and the functions it annihilates. This disconnection prevents an ``adjustment of Taylor jets" result of the type given by Bell's Lemma and thus the essential line of the Bell-Ligocka program is blocked.

To see this more explicitly, let $G=F^{-1}$. The transformation formula for the weighted Bergman projection is
$$ B_{t|z|^2}\left( JF\cdot\phi\circ F\right)=JF\cdot\left(  B_{t|G|^2}[\phi]\right)\circ F,\qquad \phi\in L^2\left(\Omega_2, e^{-t|G|^2}\right).$$
The Bergman projection $B_{t|G|^2}= B_2$ reproduces $f\in\co (\Omega_2)$. In order to find an operator $f\to Lf$ which satisfies
\begin{itemize}
\item[(i)] $B_2\left(Lf\right)=f$, and
\item[(ii)] $Lf$ vanishes to high order on $b\Omega_2$,
\end{itemize}
one is forced to consider
\begin{equation}\label{E:intro1}
e^{t|G|^2} T^p\left(e^{-t|G|^2}\cdot r^q\cdot\text{smooth}\right)
\end{equation}
where $p,q\in\Z^+$, $r$ defines $\Omega_2$, $T$ is an anti-holomorphic derivative, and smooth denotes a function in $C^\infty\left(\overline\Omega_2\right)$.
Obviously, derivatives land on $G$ in \eqref{E:intro1}---and these are the very quantities one wants to control.
The result is a vicious circle, with no boundary estimates on $F$ following from the known regularity of $B_{t|z|^2}$.

Observations of this kind seem to suggest that only estimates on the unweighted Bergman projection can be significantly connected to biholomorphic mappings. The situation changes, however, when the range of the projections are not restricted to holomorphic functions.
The purpose of this paper is to show that smoothness-to-the-boundary of $F$ can be obtained from regularity of a family of weighted, {\it non-holomorphic} projections. These projections are defined using two perturbation terms: one of them, $\tau$, shifts the space $\co (\Omega)$, and the other, $w$, weights the $L^2$ norms in the same manner as Kohn's weight mentioned above. The main result is

\begin{theorem} \label{T:main} Let $\Omega_1$ and $\Omega_2$ be smoothly bounded, pseudoconvex domains in $\C^n$ and suppose that $\Omega_1, \Omega_2$ satisfies Condition ${\fR}$. Let
 $F:\Omega_1\longrightarrow\Omega_2$ be a biholomorphic map. 

Then the components of $F$ and $F^{-1}$ extend smoothly to $\overline\Omega_1$ and
 $\overline\Omega_2$, respectively.
 \end{theorem}
 
 Condition $\fR$ contains two separate features, compatibility and regularity, about a family of twist-weight factors $(\tau, w)$ on the domains $\Omega_j$.
 It is, first of all, essential that the pairs $(\tau, w)$ be {\it Bell compatible} (see Definition \ref{D:bellcompat2}). It is also necessary that the associated family
 of twisted-weighted Bergman projections 
 $\cb^{\tau,w}_\Omega$, defined in Section \ref{S:twist}, satisfy the regularity condition given in Definition \ref{D:frakR}.

We postpone addressing the question of which domains satisfy Condition $\fR$ here and simply prove Theorem \ref{T:main},  in order to expose the twisting-weighting idea clearly. This simple idea seems to open new avenues for studying other questions in complex analysis. We hope this justifies the inclusion of some routine proofs below, e.g., Propositions \eqref{P:berg_kern}--\eqref{P:transformation}.

\medskip

Without the encouragement of several colleagues, this paper might never have appeared in manuscript form. I am especially grateful to A.-K. Herbig for
her enthusiasm about these results and gentle prodding for the past eight years to write them down. I also thank D. Varolin for insisting the results were
worthwhile and his guarantee they would not go unread. Y. Zeytuncu made a valuable observation about Section \ref{S:compat} that I happily acknowledge. And I want to thank C. L. Fefferman, J. J. Kohn, E. M. Stein, and E. J. Straube for listening to many hours of lectures about this material and for freely offering their insight on the mathematics in and behind these results.
\medskip

This paper was greatly inspired by the work in \cite{Bell81a} and \cite{BellLigocka}. Much of the proof of Theorem \ref{T:main} amounts to modest modifications of Bell's ideas.

\section{{\bf Twisted Bergman projections}}\label{S:twist}

Let $\Omega\subset\C^n$ be a bounded domain with $C^\infty$ smooth boundary, shortened to a \textit{smoothly bounded domain} below, and $r$ a smooth defining function:
$\Omega=\left\{ z\in\C^n: r(z) < 0\right\}$ and $dr\neq 0\text{ on } b\Omega=\left\{ z: r(z)=0\right\}$.
Let $\co(\Omega)$ denote the set of holomorphic functions on $\Omega$ and $C^k(\Omega)$ the $k$-times continuously differentiable functions on $\Omega$. We consider $L^2$ projections onto cosets of the form
$C^1(\Omega)\cdot\co(\Omega)$ with respect to weighted $L^2$ inner products.

To begin, if $w:\Omega\longrightarrow\R$ is a function such that $e^{-w}\in L^1_{\text{loc}}(\Omega)$, define the weighted $L^2$ function space
\begin{equation}\label{E:weightedL2}
L^2\left(\Omega, e^{-w}\right)=\left\{ f \text{ measurable on } \Omega: \int_\Omega |f|^2\, e^{-w}\, dV <\infty\right\},
\end{equation}
where $dV$ stands for the euclidean volume element. Call $w$ a weight factor.
We denote the dependence of the inner product
and norm on $w$ by a subscript: if $f, g\in L^2\left(\Omega, e^{-w}\right)$
\begin{equation}\label{E:weightedinner}
(f,g)_w=\int_\Omega f\cdot \bar g\, e^{-w}\, dV, \qquad\text{and}\qquad ||f||^2_w=(f,f)_w.
\end{equation}

Next, if $\tau:\Omega\longrightarrow\R^+$ is a positive function, belonging to $C^1(\Omega)$, define the space
\begin{equation}\label{E:twistedholo}
\co^\tau(\Omega)=\left\{ f\in C^1(\Omega):\frac{\partial}{\partial \bar z_k}\left(\sqrt\tau\cdot f\right)=0,\quad k=1,\dots, n\right\},
\end{equation}
where $\frac{\partial}{\partial\bar z_k}=\frac 12\left(\frac{\partial}{\partial x_k} +i\frac{\partial}{\partial y_k}\right)$
are the Cauchy-Riemann operators with respect to the standard coordinates $(z_1,\dots, z_n)$,
$z_k=x_k+iy_k$, on $\C^n$. Call $\tau$ a twist factor and write $\bar\partial\left(\sqrt\tau\cdot f\right)=0$ to express the vanishing of the $n$ equations in \eqref{E:twistedholo}. Clearly, $\co^\tau(\Omega)=\frac 1{\sqrt\tau}\cdot \co(\Omega)$ as sets; these are the sets we will project $L^2\left(\Omega, e^{-w}\right)$ onto. The set $\co^\tau(\Omega)$ is called the set of $\tau$-{\it twisted holomorphic functions} on $\Omega$.

\subsection{The basic inequality}\label{SS:tau_ineq}
\bigskip

For a general twist-weight pair $(\tau, w)$, let $A^2_{\tau,w}(\Omega)$ denote the $\tau$-twisted holomorphic functions in $L^2\left(\Omega, e^{-w}\right)$. If $w$ satisfies a mild integrability condition near $b\Omega$, an inequality of Bergman type holds for functions in $A^2_{\tau,w}(\Omega)$. For $\delta >0$, let $S_\delta=\left\{ z\in \Omega: |r(z)|<\delta\right\}$.

\begin{proposition}\label{T:bergineq} Let $\Omega$ be a smoothly bounded domain and $(\tau, w)$ a twist-weight pair on $\Omega$. 
Suppose that $w\in L^1_{\text loc}\left(S_\delta\right)$, for some $\delta >0$.

Then, for any compact $K\subset\Omega$, there exists a constant $C_K$ such that
$$\sup_{z\in K}|f(z)|\leq C_K\, ||f||_w,\qquad f\in A^2_{\tau,w}(\Omega).$$
\end{proposition}

\begin{remark} The constant $C_K$ also depends on the functions $\tau$ and $w$. Crucially, it is independent of
$f\in A^2_{\tau,w}(\Omega)$.
\end{remark}

\medskip

\begin{proof} If $ f\in\co^\tau(\Omega)$, then $f=\frac 1{\sqrt\tau}\cdot h$ for some $h \in\co(\Omega)$.

Choose a compact set $K_1$ such that $K\subset K_1$ and $bK_1\subset S_\delta$. Then

\begin{align*}
 \sup_{z\in K}|f(z)|&\leq \sup_{z\in K_1}\frac 1{\sqrt{\tau(z)}}\cdot\sup_{z\in K_1} |h(z)| \\
&\leq C(K_1, \tau)\cdot \sup_{z\in bK_1} |h(z)|,
\end{align*}
by the maximum principle for holomorphic functions.

Let $B(p,\eta)$ denote the euclidean ball centered at $p$ of radius $\eta$ and let $V(p,\eta)$ denote the volume of $B(p,\eta)$.
Choose $\rho >0$ such that
$B(z,\rho)\subset S_\delta$ for all $z\in bK_1$. Let $z\in bK_1$ be temporarily fixed. Since $\log |h|$ is subharmonic and $w$ is integrable on $B(z,\rho)$, we have

\begin{align*}
 2\log |h(z)|&\leq \frac 1{V(z,\rho)}\, \int_{B(z,\rho)} 2\log |h|\, dV \\
&=\int_{B(z,\rho)}\log\left(|h|^2\cdot e^{-w}\right)\frac {dV}{V(z,\rho)} +\frac 1{V(z,\rho)}\,\int_{B(z,\rho)} w\, dV.
\end{align*}
Exponentiating both sides and applying Jensen's inequality yields

\begin{align*}
 |h(z)|^2&\leq \exp\left[C(\rho, n)\int_{B(z,\rho)} w\, dV\right]
\cdot\int_{B(z,\rho)} |h|^2\, e^{-w}\,  \frac {dV}{V(z,\rho)} \\ &\leq C \left\{\sup_{x\in B(z,\rho)} \tau(x)\right\}\cdot
\int_{B(z,\rho)} |f|^2\, 
e^{-w}\, dV \\ &\leq C' ||f||^2_w.
\end{align*}
Now cover $bK_1$ by balls $B(z_1, \rho), \dots B(z_m, \rho)$. It follows that
$$ \sup_{z\in bK_1} |h(z)|\leq C''\, ||f||_w,$$
which completes the proof.

\end{proof}

\begin{remark}  The hypothesis on $w$ in Proposition \ref{T:bergineq} can be weakened.
For example, if $K_1\subset K_2\subset\dots$ are compact sets which exhaust $\Omega$ and
$\left\{ U_j\right\}$, $j\in\Z^+$, are open subsets of $\Omega$ such that $bK_j\subset U_j$, then
we need only require that $w\in L^1_{\text loc}\left(U_j\right)$, for $j\geq J$, in order to conclude
that the inequality in Proposition \ref{T:bergineq} holds. This observation shows that 
Proposition \ref{T:bergineq} holds for weights $w$ that are identically $= +\infty$ on ``rings'' accumulating to $b\Omega$ as long as there are complementary ``rings'' accumulating to the boundary where $w$ is locally integrable.
\end{remark}

\subsection{The kernel function}\label{SS:kernel}

\bigskip

From now on, consider twist-weight pairs $(\tau, w)$ with $w\in L^1_{\text loc}\left(S_\delta\right)$, for some $\delta >0$.
Proposition \ref{T:bergineq} implies that $A^2_{\tau, w}(\Omega)$ is a closed subset of $L^2\left(\Omega, e^{-w}\right)$. It also implies that for any fixed $a\in\Omega$, the evaluation functional
$$f\longrightarrow f(a),\qquad f\in A^2_{\tau, w}(\Omega)$$
is continuous in the $||\cdot||_w$ norm. The Riesz representation theorem gives, for each fixed $a\in \Omega$, a function $R_a\in A^2_{\tau, w}(\Omega)$ such that $f(a)=\left( f, R_a\right)_w$. Rewriting this, we obtain
\begin{equation}\label{E:repro}
f(a)=\int_\Omega\cb^{\tau, w}_\Omega (a,b) f(b)\, e^{-w(b)}\, dV(b), \qquad  f\in A^2_{\tau, w}(\Omega),
\end{equation}
where $\cb^{\tau, w}_\Omega (a,b)=:\overline{R_a(b)}$. This function is the $(\tau, w)$-{\it Bergman kernel} associated to $\Omega$. When the parameters $\tau, w$, and $\Omega$ are clear, we drop the super and subscripts on the kernel.

\begin{proposition}\label{P:berg_kern}The $(\tau, w)$-Bergman kernel associated to $\Omega$, 
$\cb^{\tau, w}_\Omega (a,b)=\cb(a,b)$, satisfies

\begin{itemize}
\item[(i)] \eqref{E:repro} holds
\item[(ii)] $\overline{\cb(a,\cdot)}\in A^2_{\tau, w}(\Omega)$ for each $a\in \Omega$,
\item[(iii)] $\overline{\cb(a,b)}=\cb(b,a)$.
\end{itemize}
Moreover the properties (i)-(iii) uniquely determine $\cb(a,b)$.
\end{proposition}

\begin{proof} Properties (i) and (ii) hold by definition. To see property (iii), apply \eqref{E:repro} to
$\overline{\cb(a,\cdot)}\in H^2_{\tau, w}(\Omega)$:
\begin{align*}
 \overline{\cb(a,b)}&=\int_\Omega \cb(b,s)\overline{\cb(a,s)} e^{-w(s)}\, dV(s) \\ 
&=\overline{\int_\Omega \overline{\cb(b,s)} \cb(a,s) e^{-w(s)}\, dV(s)} \\
&=\overline{\overline{\cb(b,a)}} = \cb(b,a).
\end{align*}

To verify uniqueness, suppose $K(a,b)$ is another kernel satisfying (i)-(iii). Then
\begin{align*}
\cb(a,b)&=\overline{\cb(b,a)}=\int_\Omega K(a,s)\overline{\cb(b,s)} e^{-w(s)}\, dV(s) 
\\ &=\overline{\int_\Omega K(s,a)\cb(b,s) e^{-w(s)}\, dV(s)} \\ &=\overline{K(b,a)} = K(a,b).
\end{align*}

\end{proof}

\subsection{Transformation formula}\label{SS:trans}
\bigskip

Let $\Omega_1, \Omega_2$ be smoothly bounded domains in $\C^n$ and suppose $F:\Omega_1\longrightarrow\Omega_2$ is a biholomorphic map. If $\tau$ and $w$ are functions defined on $\Omega_1$, let $\sigma =\tau\circ F^{-1}$ and $v=w\circ F^{-1}$ be the corresponding functions defined on $\Omega_2$. We want to express the relationship between the $(\tau,w)$-Bergman kernel
on $\Omega_1$ and the $(\sigma, v)$-Bergman kernel on $\Omega_2$. 

Let $JF(s)=\det \left[F'(s)\right]$ denote the determinant of the holomorphic Jacobian matrix of $F$ and use the symbol $J_{\Bbb R} F(s)$ to denote the determinant of the real Jacobian matrix of $F$, i.e., where $F$ is viewed as a diffeomorphism from ${\Bbb R}^{2n}$ to ${\Bbb R}^{2n}$. Two elementary facts are used in the proof below:
\begin{equation}\label{E:jacob}
JF\left(F^{-1}(\beta)\right)=\frac 1{JF^{-1}(\beta)},\qquad J_{\Bbb R} F^{-1}(\beta) =
JF^{-1}(\beta)\cdot\overline{JF^{-1}(\beta)}.
\end{equation}
Also, write $\cb^{*, *}_{\Omega_j} (\cdot, \cdot)=\cb^{*, *}_{j} (\cdot, \cdot)$, for $j=1,2$.

\begin{proposition}\label{P:kern_trans}
If $F:\Omega_1\longrightarrow\Omega_2$ is a biholomorphic mapping between smoothly bounded domains in $\C^n$, $(\tau,w)$ a twist-weight pair defined on $\Omega_1$, and $(\sigma,v)$ the corresponding twist-weight pair on $\Omega_2$, then
\begin{equation}\label{E:kern_trans}
\cb^{\tau, w}_{1} (a,b)=\cb^{\sigma ,v}_{2} \left(F(a),F(b)\right)\cdot JF(a)\cdot\overline{JF(b)}
\end{equation} 
for all $(a,b)\in \Omega_1\times\Omega_1$.
\end{proposition}

\begin{proof} Let $K(a,b)$ denote the function on the right hand side of \eqref{E:kern_trans}. We use Proposition \ref{P:berg_kern} to show that $K(a,b)=\cb^{\tau, w}_1(a,b)$. Property (iii) of Proposition \ref{P:berg_kern} holds for $K(a,b)$ because of the corresponding property for $\cb^{\sigma, v}_2(x,y)$. For property (ii), note first that
$$0=\bar\partial\left(\sqrt{\sigma(\cdot)}\cdot\cb^{\sigma, v}_2\left(\cdot, F(b)\right)\right) \implies
0=\bar\partial\left(\sqrt{\tau(a)}\cdot\cb^{\sigma, v}_2\left(F(a), F(b)\right)\right).$$
Since $\phi\in\co^\tau(\Omega_1)$, $\psi\in\co(\Omega_1)$ implies that $\phi\cdot\psi\in 
\co^\tau(\Omega_1)$, it follows that $K(\cdot, b)\in\co^\tau(\Omega_1)$.

For property (i), let $g\in A^2_{\tau, w}(\Omega_1)$. Applying the change of variables $b=F^{-1}(\beta)$, we have
\begin{align*}
 \int_{\Omega_1}K(a,b) g(b) e^{-w(b)}\, dV(b) = JF(a)\, \int_{\Omega_2}\cb^{\sigma, v}_2\left( F(a), \beta\right) \, &\overline{JF\left(F^{-1}(\beta)\right)}\, g\left(F^{-1}(\beta)\right)  \\
&\cdot e^{-v(\beta)}\, J_{\Bbb R}F^{-1}(\beta)\, dV(\beta).\end{align*}
However, \eqref{E:jacob} shows that the
right-hand side is
\begin{equation}\label{E:cv}
JF(a)\, \int_{\Omega_2}\cb^{\sigma, v}_2\left( F(a), \beta\right) \, \left[
 JF^{-1}(\beta)\,g\left(F^{-1}(\beta)\right) \right] 
\cdot e^{-v(\beta)}\,  dV(\beta)=M. 
\end{equation}
Since the quantity $[\dots]$ in \eqref{E:cv} is in $A^2_{\sigma, v}(\Omega_2)$ --- by the change of variables theorem --- it follows from the reproducing property of $\cb^{\sigma, v}_2$ that
\begin{align*}
 M&= JF(a)\cdot JF^{-1}\left(F(a)\right)\cdot g(a) \\ &= g(a).\end{align*}
Here \eqref{E:jacob} has been used. The uniqueness statement in Proposition \ref{P:berg_kern} now completes the proof.    
\end{proof}

The transformation formula \eqref{E:kern_trans} may also be written at the operator level. First, extend the operator in \eqref{E:repro} to all of $L^2\left(\Omega_1, e^{-w}\right)$. The $(\tau, w)$-Bergman projection is defined
$$\cb^{\tau, w}_{\Omega_1}[g](a)=\int_{\Omega_1}\cb^{\tau, w}_{\Omega_1}(a,b)g(b) e^{-w(b)}
\, dV(b),\qquad g\in L^2\left(\Omega_1, e^{-w}\right).$$
It follows from Proposition \ref{P:berg_kern} that $\cb^{\tau, w}_{\Omega_1}: L^2\left(\Omega_1, e^{-w}\right)\longrightarrow A^2_{\tau, w}(\Omega_1)$ is the orthogonal projection of $ L^2\left(\Omega_1, e^{-w}\right)$ onto $ A^2_{\tau, w}(\Omega_1)$.

\begin{proposition}\label{P:transformation} Let $F:\Omega_1\longrightarrow\Omega_2$ be biholomorphic, $(\tau, w)$ a twist-weight pair on $\Omega_1$, and 
$(\sigma ,v)$ the corresponding pair on $\Omega_2$.

Then, denoting $\cb^{\tau, w}_{\Omega_1}$ as $\cb_1$ and $\cb^{\sigma, v}_{\Omega_2}$ as $\cb_2$, 
\begin{equation}\label{E:transform}
\cb_1\left( JF\cdot\phi\circ F\right)=JF\cdot\left( \cb_2[\phi]\right)\circ F,\qquad \phi\in L^2\left(\Omega_2, e^{-v}\right).
\end{equation}
\end{proposition}

\begin{proof} Making the change of variables $\zeta= F^{-1}(\xi)$,
\begin{align*}
\cb_1\left( JF\cdot\phi\circ F\right(z))=\int_{\Omega_1}\cb_1(z,\zeta) JF(\zeta)\cdot &\phi\left(F(\zeta)\right) e^{-w(\zeta)}\, dV(\zeta) \\ = \int_{\Omega_2} 
\cb_1\left(z, F^{-1}(\xi)\right) J&F\left(F^{-1}(\xi)\right)\cdot \phi(\xi) e^{-v(\xi)} \\
&\cdot JF^{-1}(\xi)\overline{JF^{-1}(\xi)}\, dV(\xi).\end{align*}

By Proposition \ref{P:kern_trans}, this
\begin{align*}
= \int_{\Omega_2}\cb_2\left(F(z), \xi\right)&\cdot JF(z)
\cdot\overline{JF\left(F^{-1}(\xi)\right)} \\
&\cdot JF\left(F^{-1}(\xi)\right)JF^{-1}(\xi)\overline{JF^{-1}(\xi)}\, \phi(\xi) e^{-v(\xi)}\, dV(\xi)\end{align*}

\begin{align*}
=JF(z)\, \int_{\Omega_2}\cb_2&\left(F(z), \xi\right)\phi(\xi) e^{-v(\xi)} \\
&\cdot\left[\overline{JF\left(F^{-1}(\xi)\right)} JF\left(F^{-1}(\xi)\right) JF^{-1}(\xi)\overline{JF^{-1}(\xi)}\right]\, dV(\xi) \end{align*}

$$= JF(z)\cdot \left(\cb_2[\phi]\circ F\right) (z).$$
For the last equality, \eqref{E:jacob} is used to show that $[\dots]=1$. This completes the proof. 
\end{proof}

\section{{\bf The role of pseudoconvexity}}\label{S:psc}

We shall use pseudoconvexity through the following result of Diederich and Forn\ae ss, \cite{DieFor77-2}: if $\Omega$ is a smoothly bounded, pseudoconvex domain in $\C^n$, there exists a smooth defining function $\rho$ for $\Omega$ and a positive exponent $\eta$, $1\geq\eta >0$, such that $-(-\rho)^\eta$ is plurisubharmonic on $\Omega$. It is known that there is no strictly positive lower bound on $\eta$, over the class of all smoothly bounded pseudoconvex domains, for which this theorem holds; see \cite{DieFor77-2} and \cite{DieFor77-1}.

This result of Diederich-Forn\ae ss implies that, if $r$ is a defining function for a smoothly bounded pseudoconvex $\Omega$, $\log (-r)$ is quasi-invariant under biholomorphic mappings of $\Omega$. This corollary of \cite{DieFor77-2} was obtained independently by Range and Forn\ae ss:

\begin{proposition}[\cite{Fornaess79} and \cite{Range78}]\label{P:distancecomparison}
Let $\Omega_1, \Omega_2$ be smoothly bounded, pseudoconvex domains in $\C^n$ and $F:\Omega_1\longrightarrow\Omega_2$ a biholomorphic map. 

There exists an $d\in\Z^+$ such that if $r_1, r_2$ are defining functions for $\Omega_1, \Omega_2$ respectively, there are constants $C_1, C_2$ such that
\begin{equation}\label{E:distancecomparison}
C_1\cdot \left|r_1(z)\right|^d\leq \left| r_2\big(F(z)\big)\right| \leq C_2\cdot \left| r_1(z)\right|^{\frac 1d},\qquad\text{for all } z\in\Omega_1.
\end{equation}

The constants in \eqref{E:distancecomparison} depend on $F$ and the defining functions $r_1, r_2$, but are independent of $z\in\Omega_1$.
\end{proposition}

Proposition \ref{P:distancecomparison} says in particular that if $g\in C^\infty\left(\overline{\Omega_2}\right)$ vanishes on $b\Omega_2$, then $g\circ F$ must vanish on $b\Omega_1$ (though perhaps to lesser order). We will use this to estimate Sobolev norms on $\Omega_1$ by shifted Sobolev norms on $\Omega_2$. For $s\in\Z^+$, let
$$\| u\|^2_{(s)} =\sum_{|\alpha|\leq s}\int_\Omega \left|D^\alpha u\right|^2\, dV,\qquad u\in C^\infty(\Omega),$$
denote the $L^2$ Sobolev norm of order $s$. If $h\in\co(\Omega)$, this norm can be expressed using only anti-holomorphic derivatives:
\begin{equation}\label{E:holoSob}
\| h\|^2_{(s)} =\sum_{|\alpha|\leq s}\int_\Omega \left|
\frac{\partial^\alpha h}{\partial z^\alpha}\right|^2\, dV, \qquad h\in\co(\Omega).
\end{equation}
Let $W^s(\Omega)$ denote the closure of $C^\infty\left(\overline{\Omega}\right)$ in the norm $\|\cdot\|_{(s)}$ and $W^s_0(\Omega)$ denote the closure of $C^\infty_0(\Omega)$ in this norm. 

A class of multipliers on the spaces $W^s_0(\Omega)$ arises naturally.

\begin{definition}\label{D:Gtm}
Let $\Omega\subset\C^n$ be a smoothly bounded domain. For $t\in\R$, $m\in\Z^+$, and $r$ a defining function for $\Omega$, define
$$\sg_m^t(\Omega)=\left\{\mu\in C^\infty(\Omega): r^{-t+|\alpha|}\, D^\alpha\mu\in L^\infty(\Omega),\text{ for }|\alpha|\leq m\right\}.$$
Also set
$$\sg_\infty^t(\Omega)=\left\{\mu\in C^\infty(\Omega): r^{-t+|\alpha|}\, D^\alpha\mu\in L^\infty(\Omega),\forall\,\alpha\right\}.$$
Call elements in $\sg_\infty^t(\Omega)$ good $W_0^*(\Omega)$ multipliers of shift $t$.
\end{definition}

\begin{remark}
The spaces $\sg_m^t(\Omega)$ do not depend on the choice of $r$ in Definition \ref{D:Gtm}.
\end{remark}

The notation $A\lesssim B$ will henceforth express the inequality $A\leq\kappa\cdot B$ for some constant $\kappa$. The constant $\kappa$ will be independent of certain parameters, made clear in context.

\begin{lemma}\label{L:mult} Let $\mu\in\sg_m^t(\Omega)$.

\begin{itemize} 
\item[(i)] If $t\in\Z^+$ and $m\geq t$, then $\mu\in W^t_0(\Omega)$.
\item[(ii)] If $s\in\Z^+$ and $s \leq m$, 
the multiplication operator $g\longrightarrow\mu\cdot g$ maps $W_0^{s-t}(\Omega)$ to $W_0^s(\Omega)$ boundedly. 
\end{itemize}
\end{lemma}

\begin{proof} The condition $\mu\in\sg_m^t(\Omega)$ says
\begin{equation}\label{E:3.6a}
\left|D^\alpha \mu(z)\right|\lesssim |r(z)|^{t-|\alpha|},\quad |\alpha|\leq m.
\end{equation}
For (i), if $t\leq m$, \eqref{E:3.6a} implies 
\begin{equation*}
\left|D^\alpha \mu(z)\right|\lesssim |r(z)|^{t-|\alpha|},\quad |\alpha|\leq t-1.
\end{equation*}
and $D^\alpha\mu\in L^\infty(\Omega)$ if $|\alpha|=t$. Thus, $\mu\in W^t(\Omega)$ and the trace of $D^\alpha\mu$ on $b\Omega$ vanishes for $|\alpha|\leq t-1$.
Theorem 11.5 in \cite{LionsMagenes} implies that $\mu\in W_0^t(\Omega)$.

For (ii), let $g\in C^\infty\left(\overline\Omega\right)$ satisfy
\begin{equation*}
\left|D^\alpha g(z)\right|\lesssim \left\{ \begin{array}{ll} \left|r(z)\right|^{s-t-|\alpha|} & \mbox{$|\alpha|\leq s-t$} \\
1 & \mbox{$|\alpha| >s-t$}
\end{array} \right.
\end{equation*}
with constant independent of $z$ in a fixed neighborhood of $b\Omega$, i.e., $g$ vanishes to order $s-t$ on $b\Omega$.
For a fixed multi-index $\beta$, with $|\beta|\leq s-1$, it follows that
\begin{align*}
\left|D^\beta\left[\mu\cdot g\right](z)\right|&=\left|\sum_{\gamma+\delta=\beta}\frac{\beta !}{\gamma !\,\delta !}D^\gamma\mu\cdot D^\delta g\right| \\
&\lesssim \sum_{\gamma+\delta=\beta \atop |\delta|\leq s-t} |r(z)|^{t-|\gamma|}\cdot |r(z)|^{s-t-|\delta|} + \sum_{\gamma+\delta=\beta \atop |\delta| > s-t} |r(z)|^{t-|\gamma|}\\
&=|r(z)|^{s-|\beta|}.
\end{align*}
Thus $\mu\cdot g$ vanishes to order $s$ on $b\Omega$. Theorem 11.5 in \cite{LionsMagenes} implies that $\mu\cdot g\in W_0^s(\Omega)$. Since 
the set of $C^\infty\left(\overline\Omega\right)$ functions vanishing to order $s-t$ on $b\Omega$ is dense in $W^{s-t}_0(\Omega)$, the proof is complete.

\end{proof}

If $h$ is a bounded holomorphic function on $\Omega$, Cauchy's estimates imply
\begin{equation}\label{E:cauchy}
\left|D^\gamma h(z)\right|\lesssim |r(z)|^{-|\gamma|}
\end{equation}
where the constant depends only on $\Omega$, the differentiation order $\gamma$, and $\sup_\Omega |h|$. Thus $h\in\sg_\infty^0(\Omega)$. The product rule shows that $g_1\in\sg_\infty^{t_1}(\Omega), 
g_2\in\sg_\infty^{t_2}(\Omega)\implies g_1\cdot g_2\in\sg_\infty^{t_1+t_2}(\Omega)$. In particular, each component of the biholomorphic map $F=\left(f^1,\dots ,f^n\right):
\Omega_1\longrightarrow\Omega_2$ is an element of $\sg_\infty^0\left(\Omega_1\right)$ and $JF\in\sg_\infty^{-n}\left(\Omega_1\right)$.

Because of the distortion exponent $d$ in Proposition \ref{P:distancecomparison}, pullbacks of functions in $\sg_m^t\left(\Omega_2\right)$ can only be asserted to belong to shifted spaces $\sg_{\tilde m}^{\tilde t}\left(\Omega_1\right)$.

\begin{lemma}\label{L:r^p}
Let $\Omega_1, \Omega_2\subset\C^n$ be smoothly bounded, pseudoconvex domains, $F:\Omega_1\longrightarrow\Omega_2$ a biholomorphic map, and $r_2$ a defining function for $\Omega_2$. Let $d=d(F,\Omega_1,\Omega_2)$ be associated to $F$ by Proposition \ref{P:distancecomparison}.

For any $k\in\Z^+$, 
\begin{itemize}
\item[(i)] $\big( r_2\circ F\big)^{2k}\in W_0^{\frac{2k}{1+d}}(\Omega_1)$.
\item[(ii)] $\big( r_2\circ F\big)^{2k}\in \sg_k^{k/d}\left(\Omega_1\right)$.
\end{itemize}
\end{lemma}

\begin{proof} If $\beta$ is a multi-index in $\N^n$, let $D^\beta F$ denote the derivative of order $\beta$ of an (unspecified) component of $F$.
The product and chain rules imply
\begin{equation*}
D^\alpha\left(r_2^{2k}\circ F\right)=\sum_{j=1}^{|\alpha|}\sum_{\sum |\beta_j|=|\alpha|} C_{k\alpha j\beta}\big(r_2\circ F\big)^{2k-j}\cdot D^{\beta_1}F\dots D^{\beta_j}F.
\end{equation*}
for some combinatorial constants $C_*$. It follows from \eqref{E:cauchy} and Proposition \ref{P:distancecomparison}
\begin{align}\label{E:r^p1}
\left|D^\alpha\left(r_2\circ F\right)^{2k}(z)\right|&\lesssim \left| r_2\circ F\right|^{2k-|\alpha|}\cdot \left| r_1(z)\right|^{-|\alpha|} \notag \\
&\lesssim \left| r_1(z)\right|^{\frac {2k}d -|\alpha|\left(1+\frac 1d\right)}\qquad\text{for } |\alpha|\leq 2k
\end{align}
and
\begin{equation*}
\lesssim |r_1(z)|^{-|\alpha|}\qquad\text{for }|\alpha| >2k.
\end{equation*}

If $|\alpha|\leq\frac{2k}{1+d}$, this implies $\left|D^\alpha\left(r_2\circ F\right)^{2k}(z)\right|\in L^\infty(\Omega_1)$ and that $r_2\circ F$ vanishes to order at least $\frac{2k}{1+d}$, which is (i).
If $|\alpha|\leq k$, \eqref{E:r^p1} implies $\left|D^\alpha\left(r_2\circ F\right)^{2k}(z)\right|\lesssim |r_1(z)|^{\frac kd -|\alpha|}$, which is the claimed result (ii).
\end{proof}

\begin{proposition}\label{P:bdedSob}
Suppose $\Omega_1, \Omega_2$ are smoothly bounded, pseudoconvex domains in $\C^n$ and $F:\Omega_1\longrightarrow\Omega_2$ is a biholomorphic map. 
Let $\mu\in\sg_\infty^k\left(\Omega_1\right)$.

For each $s\in\Z^+$, there exists $T(s)\in\Z^+$ such that the operator $$g\longrightarrow\mu\cdot g\circ F$$ maps $W_0^{s+T(s)}\left(\Omega_2\right)$ to 
$W_0^s\left(\Omega_1\right)$ boundedly.
\end{proposition}

\begin{proof} Fix $s\in\Z^+$ and defining functions $r_1, r_2$ for $\Omega_1$ and $\Omega_2$. Let $d$ be the exponent associated to $F$ by Proposition \ref{P:distancecomparison} and let $w=F(z)$.

Let $g\in C^\infty\left(\Omega_2\right)$ vanish to order $Q$, $Q>s$ to be determined, on $b\Omega_2$:
\begin{equation}\label{E:bS1}
\left| D^\delta g(w)\right|\lesssim \left|r_2(w)\right|^{Q-|\delta|}, \quad |\delta|\leq Q,\qquad w\text{ near }b\Omega_2.
\end{equation}

If $F=\left(f^1,\dots ,f^n\right)$, each $f^k$ satisfies \eqref{E:cauchy}. The chain rule gives, for $|alpha|\leq Q$,
\begin{equation*}
\left| D^\alpha\left[g\circ f\right](z)\right|\lesssim\sum_{\beta: |\beta|=|\alpha|}\left|D^\beta g\big(F(z)\big)\right|\cdot a_1\cdot\dots\cdot a_{|\beta|},
\end{equation*}
where each factor $a_l\in\sg_\infty^{-l}\left(\Omega_1\right), l=1,\dots ,|\beta|$, by \eqref{E:cauchy}. Combining this with \eqref{E:bS1} and using \eqref{E:distancecomparison} yields
\begin{align*}
\left| D^\alpha\left[g\circ F\right](z)\right|&\lesssim \left|r_2(w)\right|^{Q-|\alpha|}\cdot \left|r_1(z)\right|^{-|\alpha|} \\
&\lesssim \left| r_1(z)\right|^{\frac Qd-\left(1+\frac 1d\right)|\alpha|}.
\end{align*}
Leibniz's rule thus implies
\begin{align}\label{E:bS2}
\left|D^\beta\left[\mu\cdot g\circ F\right](z)\right|&=\left|\sum_{\gamma+\delta=\beta}\frac{\beta !}{\gamma !\,\delta !}D^\gamma\mu (z)\cdot
 D^\delta\left[ g\circ F\right](z)\right| \notag\\
 &\lesssim \sum_{\gamma+\delta=\beta}\left| r_1(z)\right|^{-k-|\gamma|}\cdot\left| r_1(z)\right|^{\frac Qd -|\delta|-\frac{|\delta|}d}.
\end{align}

If $Q\geq dk +(d+1)s+1$, it follows from \eqref{E:bS2} that $\mu\cdot g\circ F$ vanishes to order $\geq s$ on $b\Omega_1$ and, consequently, 
$\mu\cdot g\circ F\in W_0^s\left(\Omega_1\right)$ by Theorem 11.5 of \cite{LionsMagenes}. If $T(s)\geq dk+(d-2)s+1$, the conclusion follows from the fact that functions in $C^\infty\left(\Omega_2\right)$ satisfying \eqref{E:bS1} are dense in $W_0^Q\left(\Omega_2\right)$.

\end{proof}

\section{{\bf A variation on Bell's operator}}\label{S:bell}

Constructions in this section occur on a single domain, so notation is simplified. Let $\Omega$ be a smoothly bounded domain in $\C^n$, given by a smooth defining function $r$.
The twist factor will be denoted by $T$ and the weight factor by $W$: $f\in\co^T(\Omega)$ if $\bar\partial\left(\sqrt T\cdot f\right) =0$ and $f\in L^2\left(\Omega, e^{-W}\right)$ if $||f||_W<\infty$ where $||\cdot||_W$ is defined in \eqref{E:weightedinner}. The $(T,W)$-Bergman projection on $\Omega$, $\cb^{T,W}_\Omega$, is denoted $\cb$.

Several spaces of functions that ``vanish on $b\Omega$'' arise in the analysis. For $L\in\Z^+$, let
$$V^L(\Omega)=\left\{ g: g=b\cdot r^L, \text{for some } b\in L^\infty(\Omega)\right\}$$ 
denote the bounded functions that vanish to order $L$ on $b\Omega$.  For $\ell\in{\Bbb Z}^+\cup\{0\}$, let
$$\sB_\ell(\Omega)=\big\{b\in C^\infty(\Omega): D^\alpha b\in L^\infty(\Omega),\quad |\alpha|\leq\ell\big\}.$$
Then, for $L\in\Z^+$, define subsets of $V^L$ by
$$\cv_\ell^L(\Omega)=\big\{ g: g=b\cdot r^L,\, b\in\sB_\ell(\Omega)\big\}.$$
The spaces $\cv_\ell^L(\Omega)$ are sometimes written $\cv_\ell^L$ below, when $\Omega$ is fixed.
Clearly, $V^L=\cv^L_0\supset\cv^L_1\supset\dots$.

Define the differential operators 
$$D_k= e^W\, \sqrt T\, \frac\partial{\partial\bar z_k},\qquad k=1,\dots, n.$$ 
The initial observation is that the image of $\cv^1_1$ under the operators $D_k$
is orthogonal to $\overline{\co^T(\Omega)}$ in $L^2\left(\Omega, e^{-W}\right)$.

\begin{lemma}\label{L:var1} Let $\Omega$ be a smoothly bounded domain and $(T,W)$ a twist-weight pair on $\Omega$ with $T\in L^\infty(\Omega)$. 

If $h\in\overline{\co^T(\Omega)}$ and $f\in\cv^1_1(\Omega)$, then
\begin{equation}\label{E:ortho1}
\big(h, D_k(f)\big)_W=0\qquad\text{for } k=1,\dots, n.
\end{equation}
\end{lemma}

\begin{proof} Let $h\in\overline{\co^T(\Omega)}$, i.e., $\partial\left(\sqrt{T}h\right)=0$. Assume, temporarily, that $h\in L^\infty(\Omega)$ as well. If $f\in\cv^1_1$, integration by parts yields
\begin{align}\label{E:ortho2}
\big( h,D_k(f)\big)_W &=\int_\Omega h\cdot\left(e^W\sqrt{T}\frac{\partial\bar f}{\partial z_k}\right)\, e^{-W}\, dV \notag \\
&=\, -\int_\Omega\frac{\partial}{\partial z_k}\left(\sqrt{T}h\right)\cdot\bar f =0.
\end{align}
No boundary term occurs because $f\in\cv^1_1$ and $\sqrt{T}h\in L^\infty(\Omega)$; the integrand in \eqref{E:ortho2} vanishes since $h\in\overline{\co^T(\Omega)}$.

For general $h\in\overline{\co^T(\Omega)}$, a limiting argument is used. A partition of unity reduces the problem to showing
\eqref{E:ortho1} holds for $f\in V$ supported near some $p\in b\Omega$. Let $\nu$ denote the outward unit normal to $b\Omega$ at $p$. 
For $h\in\overline{\co^T(\Omega)}$ given, set $h_\epsilon(z) = h(z-\epsilon \,\nu)$. If $\epsilon>0$ is small, $h_\epsilon$ is well-defined and belongs to 
$\overline{\co^T(\Omega)}\cap L^\infty(\Omega)$. If the support of $f$ is sufficiently small, the above integration by parts argument applies, giving $\big(h_\epsilon, D_kf\big)_W=0$.
However since $\sqrt{T},\frac{\partial\bar f}{\partial z_k}\in L^\infty(\Omega)$,

$$\big(h, D_k(f)\big)_W=\lim_{\epsilon\to 0}\big(h_\epsilon,  D_k(f)\big)_W$$
by dominated convergence. Thus \eqref{E:ortho1} holds under the stated hypothesis.
\end{proof}

A reformulation in terms of the null space of the $(T,W)$-Bergman projection on $\Omega$ is convenient:

\begin{corollary}\label{C:var2} Under the hypothesis of Lemma \ref{L:var1},

\begin{equation}\label{E:ortho3}
\cb\left[ \bar D_k f\right] =0,\quad k=1,\dots, n,\qquad\text{for any }f\in\cv^1_1(\Omega).
\end{equation}
\end{corollary}

\begin{proof} Recall that the $(T,W)$-Bergman kernel satisfies $\cb(a,\cdot)\in\overline{\co^T(\Omega)}$ and 
$$\cb\left[ \bar D_k f\right](a) =\left(\cb (a,\cdot), D_k\bar f\right)_W.$$
\eqref{E:ortho3} now follows from Lemma \ref{L:var1}.
\end{proof}

The main result of this section generalizes Lemma 2 in \cite{Bell81a} (see, also, \cite{Bell79} Lemma 3.1).

\begin{proposition}\label{P:mainbell} Let $\Omega$ be a smoothly bounded domain and $(T,W)$ a twist-weight pair on $\Omega$ with $T\in L^\infty(\Omega)$. 
\smallskip

Let $M\in\Z^+$ and consider functions $\phi$ of the form
$$\phi=e^W\sqrt T\cdot \left(b\cdot s\right), \qquad\text{for } b\in\sB_M(\Omega)\text{ and } s\in C^\infty(\overline\Omega).$$
For any $m\in\Z^+$, $m\leq M$,
there exists a function $\psi^m$ such that

\begin{itemize}
\item[(i)] $\cb\left[\phi\right]= \cb\left[\psi^m\right],$ and
\smallskip
\item[(ii)] $e^{-W}\frac 1{\sqrt T}\cdot \psi^m\in\cv_{M-m}^m$.
\end{itemize}
\end{proposition}

\begin{proof} 

The analysis occurs near $b\Omega$, where the defining function $r$ can be used as a coordinate. The function $s$, the ``smooth part'' of $\phi$, then locally has a Taylor expansion in powers of $r$ up to $b\Omega$. $\psi^m$ is constructed by removing terms of order $< m$ from this expansion; Corollary \ref{C:var2} implies the removed terms are in 
$\text{Null}(\cb)$.

For any $p\in b\Omega$, there exists a neighborhood $U$ and at least one anti-holomorphic derivative, say $\frac{\partial r}{\partial z_{(p)}}= r_{z_{(p)}}$, that is  
non-vanishing on $U$, since $dr\neq 0$ on $b\Omega$. Fix $p\in b\Omega$ and such a neighborhood $U$. Initially, suppose the smooth part of $\phi$ is supported in $U$,
 i.e., $\phi =e^{W}\sqrt T\cdot\big( b\cdot s\big)$ with $s\in C^\infty(\overline\Omega)\cap C^\infty_0(U)$. The function
$\psi^m$ will be of the form
\begin{equation}\label{E:bell1}
\psi^m=\phi -e^W\sqrt T\cdot\frac\partial{\partial z_{(p)}}\left(\sum_{j=1}^m g_j\cdot r^j\right),
\end{equation}
for functions $g_1,\dots, g_m$ to be chosen. Corollary \ref{C:var2} says that $\cb\left[\psi^m\right]=\cb\left[\phi\right]$, for any choices of $g_l\in\sB_1(\Omega)$. The goal is to chose $g_j$ such that $e^{-W}\frac 1{\sqrt T}\cdot \psi^m\in\cv^m_{M-m}$; the proposed form \eqref{E:bell1} reduces this to a two-term recursion relation. Choose $g_1$ satisfying
\begin{equation}\label{E:bell2}
b\cdot s-g_1\cdot r_{z_{(p)}} =0
\end{equation}
to force the $r^0$ term in \eqref{E:bell1} to vanish
(throughout, $r^t$ denotes $r$ raised to the $t$-power). It then follows that
\begin{equation*}
g_1(z)=\begin{cases}\frac{b\cdot s}{r_{z_{(p)}}}(z) &\text{ if $z\in U$} \\
0 & \text{if $z\notin U$}
\end{cases}
\end{equation*}
determines $g_1$ as an element in $\sB_M(\Omega)$, since $r_{z_{(p)}}\neq 0$ in $U$, $s$ is supported in $U$, and $b\in\sB_M(\Omega)$.

Next, choose $g_2$ satisfying
$$\frac{\partial g_1}{\partial z_{(p)}}+ 2g_2\cdot r_{z_{(p)}}=0,$$
in order to annihilate the $r^1$ term in \eqref{E:bell1}. This equation says $g_2\sim \frac{\partial g_1}{\partial z_{(p)}}\sim \frac{\partial b}{\partial z_{(p)}}$ in $U$, where $\sim$ denotes equality up to factors in $C^\infty\left(\overline{\Omega}\right)$.
As before, setting $g_2=0$ outside $U$ determines $g_2$ as an element of $\sB_{M-1}(\Omega)$.
Continue by choosing $g_j$, $j=3,\dots m$, so that
$$\frac{\partial g_{j-1}}{\partial z_{(p)}}+ j\, g_j\cdot r_{z_{(p)}}=0.$$
The result is that
\begin{equation}\label{E:bell3}
\psi^m=-e^W\sqrt T\cdot\left(\frac{\partial g_m}{\partial z_{(p)}}\right)\cdot r^m,
\end{equation}
with
\begin{equation*}
\frac{\partial g_m}{\partial z_{(p)}}\sim \frac{\partial^m b}{\partial z_{(p)}^m}\in\sB_{M-m}(\Omega).
\end{equation*}

To globalize this, let $U_1,\dots, U_N$ be neighborhoods of $p_1,\dots, p_N\in b\Omega$ such that

\begin{itemize}
\item[(i)] $r_{z_{(k)}}\neq 0$ on $U_k$, and
\item[(ii)] $b\Omega\subset\cup_{\ell =1}^N U_l$.
\end{itemize}
Let $\left\{\chi_\ell\right\}_{\ell =1}^N$ be a partition of unity subordinate to $\left\{U_\ell\right\}_{\ell =1}^N$.
On each $U_k$, the above construction yields functions $g_1^{(k)},\dots, g_m^{(k)}$, $g_\ell\in\sB_{M-\ell+1}(\Omega)$ such that
$\chi_k\phi-e^W\sqrt T\cdot\frac{\partial}{\partial z_{(k)}}\left(\sum_{\ell=1}^m g_\ell^{(k)}\cdot r^\ell\right)\,\,\in e^W\sqrt T\cdot\cv^m_{M-m}.$
 Thus
$$\psi^m =\phi -e^W\sqrt T\cdot\sum_{k=1}^N
\frac{\partial}{\partial z_{(k)}}\left(\sum_{\ell=1}^m g_\ell^{(k)}\cdot r^\ell\right)$$
gives the desired function. 
\end{proof}

\section{{\bf Compatibility of twist-weight pair}}\label{S:compat}

We first introduce some auxiliary classes of functions:

\begin{definition}
Let $\Omega$ be a smoothly bounded domain. For $t\in\R$ and $m\in\Z^+$ define
$$\sa_m^t(\Omega)=\left\{f\in C^\infty(\Omega): f\in L^\infty(\Omega),\text{ and } \left|D^\alpha f(z)\right|\lesssim |r(z)|^t, 0<|\alpha|\leq 2m\right\}.$$

If $f\in\sa_m^t$ for all $t\in\R^+$, write $f\in\sa_m^\infty(\Omega)$.
\end{definition}

\begin{remark}
Clearly, $\{\text{constants}\}\subset \sa_m^\infty(\Omega)$ for all $m\in\Z^+$. However, $\sa_m^\infty(\Omega)$ contains other functions, e.g., any function of the form $s\cdot e^{\frac 1r}$ for $s\in C^\infty\left(\overline\Omega\right)$.

If $\Omega$ is pseudoconvex, the classes $\sa_m^t(\Omega)$ are quasi-invariant under biholomorphic maps. If $F:\Omega_1\longrightarrow\Omega_2$ is a biholomorphism between smoothly bounded, pseudoconvex domains $\Omega_1, \Omega_2$ and $f\in \sa_m^t(\Omega_1)$, then $f\circ F^{-1}\in\sa^m_{\tilde t}(\Omega_2)$ for some $\tilde t$.\footnote{This argument gives $\tilde t << t$, since the estimate \eqref{E:cauchy} is used. After Theorem \ref{T:main} is known, it follows that the classes $\sa_m^t(\Omega)$ are biholomorphic invariants.} In particular, if $f\in \sa_m^t(\Omega_1)$ 
\textit{for large enough} $t$ (depending on $F$), then $f\circ F^{-1}\in \sa_m^0(\Omega_2)$.
\end{remark}

Return to the set-up in Section \ref{S:twist}: $\Omega_1, \Omega_2$ are smoothly bounded domains, 
$F:\Omega_1\longrightarrow\Omega_2$ is a biholomorphic map, $(\tau,w)$ a twist-weight pair on $\Omega_1$, and $\sigma=\tau\circ F^{-1}$, $v=w\circ F^{-1}$ the corresponding pair on $\Omega_2$. Denote the twisted-weighted Bergman projections $\cb^{\tau, w}_{\Omega_1}$, $\cb^{\sigma, v}_{\Omega_2}$ by $\cb_1$ and $\cb_2$.

The initial goal is to formulate conditions on the pair $(\tau, w)$ on $\Omega_1$ that ensure $\frac{e^{-v}}{\sigma}\in \sB_*(\Omega_2)$.

\begin{definition}\label{D:bellcompat1}
Let $\Omega\subset \C^n$ be a smoothly bounded domain, $t\in\R$, and $m\in\Z^+\cup \{0\}$.
A  twist-weight pair $(\tau, w)$ is called $\sB$-compatible to index $(m,t)$  if $\frac{e^{-w}}{\tau}\in\sa_m^t(\Omega)$.
\end{definition}

\begin{lemma}\label{L:bellcompat} Let $\Omega_1, \Omega_2\subset\C^n$ be smoothly bounded, pseudoconvex domains, $F:\Omega_1\longrightarrow\Omega_2$ a biholomorphic map, and $d=d(F,\Omega_1,\Omega_2)$ be the distortion exponent given by Proposition \ref{P:distancecomparison}.

For $L\in\Z^+$, if $(\tau,w)$ is $\sB$-compatible to index $(L,2Ld)$ on $\Omega_1$, then $(\sigma,v)$ is $\sB$-compatible to index $(L,0)$ on $\Omega_2$. In particular,
$\frac{e^{-v}}\sigma\in\sB_{2L}(\Omega_2)$ if $(\tau,w)$ is $\sB$-compatible to index $(L,2Ld)$ on $\Omega_1$.
\end{lemma}

\begin{proof}
This follows directly from \eqref{E:cauchy}, Proposition \ref{P:distancecomparison}, and the definition of the spaces $\sa_m^t(\Omega)$.
\end{proof}

The condition of $\sB$-compatibility connects the type of functions reproduced by $\cb_j$, $j=1,2$, and the multiplier in front of $b\cdot s$ in Proposition \ref{P:mainbell}. This gives

\begin{proposition}\label{P:mainstep}
Let $\Omega_1, \Omega_2\subset\C^n$ be smoothly bounded, pseudoconvex domains and $F:\Omega_1\longrightarrow\Omega_2$ a biholomorphic map with distortion exponent $d=d(F,\Omega_1,\Omega_2)$. 

Let $t\in\Z^+$ be given and suppose
$(\tau, w)$ are $\sB$-compatible to index $(L,2Ld)$, with $L=(2n+t)(1+d)$, on $\Omega_1$.

For any $g\in A^\infty\left(\Omega_2\right)$, there exists a function $G\in W^t_0\left(\Omega_1\right)$ such that

\begin{equation}\label{E:mainstep}
JF\cdot g\circ F=\sqrt\tau\,\cb_1\left[e^w\sqrt\tau\,\cdot G\right].
\end{equation}
\end{proposition}

\begin{proof} For given $g\in A^\infty\left(\Omega_2\right)$, set $\phi=\frac 1{\sqrt\sigma}\cdot g$. Note that $\phi\in A^2_{\sigma, v}\left(\Omega_2\right)$ 
since $(\sigma, v)$ are $\sB$-compatible to index $(*,0)$.
Rewrite 
$\phi$ as
\begin{equation}\label{E:naturalB_L}
\phi= e^v\sqrt\sigma\,\left(\frac{e^{-v}}\sigma\cdot g\right).
\end{equation}
Since $(\tau, w)$ are $\sB$-compatible to index $(L,2Ld)$, it follows from Lemma \ref{L:bellcompat} that $\frac{e^{-v}}\sigma\in\sB_{2L}(\Omega_2)$.
Proposition \ref{P:mainbell} applied for $M=2L, m=L$ gives $\psi^L$ such that $\cb_2\left[\psi^L\right]=\cb_2[\phi]$ and
\begin{equation}\label{E:mainstep1}
e^{-v}\frac 1{\sqrt\sigma}\,\cdot \psi^L\in\cv^{L}_{L}\left(\Omega_2\right).
\end{equation}
Since $\cb_2$ reproduces $\phi$, \eqref{E:transform} gives
\begin{align*}
JF\cdot\frac 1{\sqrt\tau}\cdot\left[ g\circ F\right]=&JF\cdot\left( \cb_2[\phi]\circ F\right) \\
=& JF\cdot\left( \cb_2\left[\psi^L\right]\circ F\right) \\
=&\cb_1\left( JF\cdot\left(\psi^L\circ F\right)\right).
\end{align*}

\eqref{E:mainstep1} says that $\psi^L=e^v\sqrt\sigma\left(b\cdot r_2^L\right)$ where $b\in\sB_L\left(\Omega_2\right)$. Inserting this above and simplifying, we obtain
\begin{equation*}
JF\cdot g\circ F=\sqrt\tau\,\cb_1\left[e^w\sqrt\tau\,\cdot G\right],
\end{equation*}
where $G=JF\cdot (b\circ F)\cdot (r_2^L\circ F)=:g_1\cdot g_2\cdot g_3$. It remains to verify that $G\in W_0^{t}\left(\Omega_1\right)$. However, $g_3\in W_0^{\frac L{1+d}}\left(\Omega_1\right)$ by Lemma \ref{L:r^p}. Also, 
$g_1\in\sg_\infty^{-n}\left(\Omega_1\right)$ and $g_2\in\sg_L^{-n}(\Omega_1)$, as noted above Lemma \ref{L:r^p}. It follows from Lemma \ref{L:mult} that
$G\in W_0^{t}\left(\Omega_1\right)$ as claimed.
\end{proof}

We shall also require that the multiplier in front of $G$, on the right-hand side of \eqref{E:mainstep}, be a good $W^*_0(\Omega_1)$ multiplier.
This is a second, separate notion of compatibility on $(\tau,w)$. Both compatibility notions are combined in the next definition, formulated on a family of twist-weight pairs in order that the conclusion of Proposition \ref{P:mainstep} holds as the biholomorphism $F$ varies. Note that the second compatibility condition is required to hold uniformly in the family.

\begin{definition}\label{D:bellcompat2}
Let $\Omega\subset \C^n$ be a smoothly bounded domain. A family of twist-weight pairs
$\cf=\left\{(\tau_i, w_i), i\in\Z^+\right\}$, on $\Omega$ is Bell compatible if

\begin{itemize}
\item[(i)] for any $(A,B)\in\Z^+\times\Z^+$, there exists $j(A,B)$ such that $(\tau_j,w_j)\in\cf$ are $\sB$-compatible to index $(A,B)$ when $j\geq j(A,B)$,
\item[(ii)] there exists $K$ such that $e^{w_i}\,\sqrt{\tau_i}\in\sg_\infty^K(\Omega)$ for all $i\in\Z^+$.
\end{itemize}
\end{definition}

\begin{corollary}\label{C:mainstep}
Let $\Omega_1$ be a smoothly bounded pseudoconvex domain in $\C^n$ and $\cf=\left\{(\tau_i,w_i), i\in\Z^+\right\}$ be a Bell compatible family on $\Omega_1$.

Let $s\in\Z^+$. Then for any biholomorphic map $F:\Omega_1\longrightarrow\Omega_2$ onto a smoothly bounded pseudoconvex domain $\Omega_2$, there exists $j(s)\in\Z^+$ such that: for any $h\in A^\infty\left(\Omega_2\right)$, there exists a function $H\in W^s_0\left(\Omega_1\right)$ such that

\begin{equation}\label{E:mainstep2}
JF\cdot h\circ F=\sqrt\tau\,\cb^{\tau_j,w_j}_1\left[H\right],
\end{equation}
for any $(\tau_j,w_j)\in\cf$ with $j\geq j(s)$.
\end{corollary}

\begin{proof} Let $K$ be the constant given by Definition \ref{D:bellcompat2} (i) for the family $\cf$.
By Proposition \ref{P:mainstep}, any twist-weight pair $(\tau,w)$ that are $\sB$-compatible to sufficiently high order
cause \eqref{E:mainstep} to hold $G\in W^{s+K}_0(\Omega)$.  
Set $H=e^w\sqrt\tau\cdot G$.
Lemma \ref{L:mult} implies that $H$ belongs to $W^{s}_0(\Omega)$, which completes the proof.

\end{proof}

The trivial family $\cf=\left\{(\tau_i,w_i)=(1,0), \text{ for all }i\in\Z^+\right\}$ is obviously Bell compatible, since the constant function $1=e^0\sqrt{1}$ belongs $\sg_\infty^K$ for any $K\in\R$ and $1=\frac{e^{-0}}{1}$ is $\sB$-compatible to any index in $\Z^+\times\Z^+$. Other Bell compatible families may be obtained, e.g., by taking a fixed function $\tau$ such that $\frac 1{\sqrt\tau}\in\sg_\infty^0(\Omega)$ and setting $w_i=-\log\left(\tau\left(1+(-r)^i\right)\right)$, $i\in\Z^+$.

\section{{\bf Proof of main theorem}}\label{S:proof}

The regularity hypothesis in Theorem \ref{T:main} can now be stated:

\begin{definition}\label{D:frakR} A smoothly bounded domain $\Omega$ satisfies \textit{Condition} $\fR$ if there exists a Bell compatible family of twist-weight pairs $\left\{(\tau_j,w_j),\,\, j\in\Z^+\right\}=\cf$ satisfying the following:

for each $s\in\Z^+$, there exist $m(s), \ell(s)\in\Z^+$ such that for all $m\geq m(s)$

\begin{equation}\label{E:frakR_reg}
\sqrt{\tau_m}\cdot\cb^m:W^{s+\ell(s)}_0(\Omega)\longrightarrow W^s(\Omega),
\end{equation}
where $\cb^m=\cb^{\tau_m,w_m}_\Omega$.
\end{definition}

Definition \ref{D:frakR}  reduces  to Condition R of Bell-Ligocka when $\cf$ is taken to be the single pair $(1,0)$.

\begin{proof}[Proof of Theorem \ref{T:main}] Let $F:\Omega_1\longrightarrow\Omega_2$ be a fixed biholomorphism and $\cf=\left\{(\tau_j,w_j)\right\}$ the Bell compatible
family on $\Omega_1$. 

Let $I\in\Z^+$ be arbitrary. First apply Corollary \ref{C:mainstep} for $s=I+\ell(I)$ where $\ell(I)$ is given by Condition $\fR$.
Set $N=\max(m(I),j(I))$, $m(I)$ given by Condition $\fR$ and $j(I)$ given by Corollary \ref{C:mainstep}. Take any $(\tau_j,w_j)\in\cf$ with $j\geq N$. It follows from  \eqref{E:frakR_reg} that

\begin{equation}\label{E:proof1}
\left\|JF\cdot h\circ F\right\|_{(I)} <\infty
\end{equation}
for any $h\in A^\infty(\Omega_2)$. Since $I$ was arbitrary, Sobolev's lemma implies $JF\cdot h\circ F\in C^\infty\left(\overline{\Omega}_1\right)$. Choosing $h\equiv 1$ yields $JF\in C^\infty\left(\overline{\Omega}_1\right)$.

The same argument on $F^{-1}:\Omega_2\longrightarrow\Omega$, using the Bell compatible family $\tilde\cf =\left\{(\tilde\tau_j,\tilde w_j)\right\}$ on $\Omega_2$, shows that $JF^{-1}\in C^\infty\left(\overline{\Omega}_2\right)$. Note that $\tilde\cf$ is not necessarily the family $\left\{(\sigma_j,v_j)\right\}$ with $(\sigma,v)$ associated to $(\tau, w)$ as previously. Since $JF^{-1}$ is smooth up to $b\Omega_2$, it follows that $JF\neq 0$ on $\overline{\Omega}_1$.

Apply \eqref{E:proof1} to the coordinate functions, $h_k(w_1,\dots , w_n)=w_k, k=1,\dots , n$, to obtain $JF\cdot f_k \in C^\infty\left(\overline{\Omega}_1\right)$,
where $F=\left(f_1,\dots , f_n\right)$. Dividing out the non-vanishing factor $JF$ yields $f_k \in C^\infty\left(\overline{\Omega}_1\right)$, $k=1,\dots , n$. Since the argument is reversible, we also obtain $\left(F^{-1}\right)_k \in C^\infty\left(\overline{\Omega}_2\right)$, which completes the proof.

\end{proof}

 \vskip .5cm
\bibliographystyle{plain}
\bibliography{mcneal}  
\end{document}